\documentclass[12pt]{article}

 \usepackage{a4}
 \usepackage{amsmath}
 \usepackage{amssymb}
\newtheorem{thm}{Theorem}[section]

 \newtheorem{prop}[thm]{Proposition}

 \newtheorem{lemma}[thm]{Lemma}

\DeclareMathOperator{\diam}{diam}
\newcommand{\Int}{\int_0^1}

\title{On a fixed point in the metric space of normalized Hausdorff moment sequences}
 \author{Christian Berg\footnote{Corresponding author}\; and 
Maryam Beygmohammadi} 
 \date{\today}

\begin{document}

 \maketitle 

\begin{abstract} We show that the transformation $(x_n)_{n\ge 1}\to
  (1/(1+x_1+\ldots+x_n))_{n\ge 1}$ of the compact set of sequences
  $(x_n)_{n\ge 1}$ 
of numbers from the unit interval $[0,1]$ has a unique fixed point, which is attractive. The fixed point
  turns out to be a Hausdorff moment sequence studied in \cite{B:D2}.
\end{abstract}


    
\section{Introduction}

Let $\mathcal K=[0,1]^{\mathbb N}$ denote the product space of
sequences $(x_n)=(x_n)_{n\ge 1}$ of numbers from the unit interval $[0,1]$.

We consider a transformation $T:\mathcal K\to \mathcal K$ defined by
\begin{equation}\label{eq:trans}
(T(x_n))_n=\frac{1}{1+x_1+\ldots+x_n},\quad n\ge 1.
\end{equation}
Since $T$ is a continuous transformation of the compact convex set
$\mathcal K$
in the space $\mathbb R^{\mathbb N}$ of real sequences
equipped with the product topology, it
has a fixed point $(m_n)$ by Tychonoff's extension of Brouwer's fixed
point theorem. Furthermore, it is clear by \eqref{eq:trans} that the
fixed point $(m_n)$ is uniquely determined by the equations
\begin{equation}\label{eq:fix}
(1+m_1+\ldots +m_n)m_n=1,\quad n\ge 1.
\end{equation}                                    
Therefore
\begin{equation}\label{eq:fixrec}
m_{n+1}^2+\frac{m_{n+1}}{m_n}-1=0,
\end{equation}
giving
$$
m_1=\frac{-1+\sqrt 5}{2},\quad m_2=\frac{\sqrt{22+2\sqrt 5}-\sqrt
5-1}{4}, \ldots\,.
$$  
 Berg and Dur{\'a}n studied this fixed point in \cite{B:D2}, and it was
proved that $m_0=1,m_1,m_2,\ldots$ is a normalized Hausdorff moment
sequence, i.e., of the form 
\begin{equation}\label{eq:Ha}
m_n=\Int x^n\,d\tau(x),
\end{equation}
where $\tau$ is a probability measure on the interval $[0,1]$.
 For details about the measure $\tau$, see \cite{B:D2} and \cite{B:D3}. We
mention that $\tau$ has an increasing and convex  density with respect to Lebesgue measure. 

There is no reason a priori that the fixed point $(m_n)$ should be a
Hausdorff moment sequence, but the motivation for the study of $T$
came from the theory of moment sequences because of the following
theorem from \cite{B:D1}:

\begin{thm} \label{thm:HtoH}
Let $(a_n)_{n\ge 0}$ be a Hausdorff moment sequence of a measure $\mu
\not =0$ on $[0,1]$. Then the sequence
$(b_n)_{n\ge 0}$ defined by $b_n=1/(a_0+\ldots +a_n)$ is again a
Hausdorff moment sequence, and its associated measure $\nu =\widehat{T}(\mu
)$ has the properties $\nu (\{ 0\})=0$ and
\begin{equation}\label{eq:assomeas}
\Int\frac{1-t^{z+1}}{1-t}\,d\mu (t)\Int t^z\,d\nu (t)=1\quad
\mbox{for}\quad \Re z\ge 0.
\end{equation}
\end{thm}

Let $\mathcal H$ denote the set of normalized Hausdorff moment
sequences, where we throw away the zero'th moment which is always 1, i.e.,
\begin{eqnarray}\label{eq:H}
\mathcal H=\{(a_n)_{n\ge 1}\;\mid\; a_n=\int_0^1 x^n\,d\mu(x),\;\mu([0,1])=1\}.
\end{eqnarray}
Clearly, $\mathcal H\subset \mathcal K$ and by Theorem~\ref{thm:HtoH}
we have $T(\mathcal H)\subseteq\mathcal H$. It is easy to see that $\mathcal
H$ is a compact convex subset of $\mathcal K$, e.g., by using
Hausdorff's 1921 characterization of Hausdorff moment sequences as
completely monotonic sequences, i.e., sequences $(a_n)_{n\ge 0}$ satisfying
\begin{equation}\label{eq:Hausdorff}
\sum_{k=0}^m (-1)^k \binom mk a_{n+k}\ge 0 \quad\mbox{for}\quad
m,n\geq 0.
\end{equation}
See \cite{Ak} for details.

It was proved in Theorem 2.3 in \cite{B:D2} that $(m_n)$ is an attractive
fixed point of the restriction of $T$ to $\mathcal H$. This proof was
a direct proof not building on any classical results on attracting
fixed points. For classical fixed point theory see \cite{S}. 

The purpose of the present paper is to prove the following extension
of this:

\begin{thm}\label{thm:attrac} The unique fixed point $(m_n)$ of the
  transformation $T:\mathcal K\to \mathcal K$ given by
  \eqref{eq:trans} is attractive.
\end{thm}

\section{Proofs and complements}

The product topology  on the vector space $\mathbb
R^{\mathbb N}$ is induced by the metric
\begin{eqnarray*}
d((a_n),(b_n))=\sum_{n=1}^\infty 2^{-n}\min\{|a_n-b_n|,1\}\quad
\mbox{for}\quad (a_n),(b_n)\in\mathbb R^{\mathbb N},
\end{eqnarray*}
which makes it a Fr\'echet space.

On the compact subset $\mathcal K=[0,1]^{\mathbb N}$ the expression for
the metric is simplified to
\begin{eqnarray}\label{eq:metricK}
d((a_n),(b_n))=\sum_{n=1}^\infty 2^{-n}|a_n-b_n|\quad
\mbox{for}\quad (a_n),(b_n)\in\mathcal K.
\end{eqnarray}

Before we find the best Lipschitz constant for $T$, let us introduce
some notation.
For $0\le a\le 1$ let $\underline{a}=(a^n)_{n\ge 1}$, so
$\underline{0}=0,0,\ldots$ and $\underline{1}=1,1,\ldots$. Clearly,
$T(\underline{0})=\underline{1}$, $T(\underline{1})=1/2,1/3,\ldots$,
while
$$
\left(T\left(\frac{1}{n+1}\right)\right)_n=\frac{1}{H_{n+1}},\quad\mbox{where}\quad H_n=\sum_{k=1}^n\frac{1}{k}.
$$
The numbers $H_n$ are called the harmonic numbers. In
\cite{B:D1} it is proved that
$$
\frac{1}{H_{n+1}}=\Int x^n\left(\sum_{p=0}^\infty \alpha_p
  x^{-\xi_p}\right)\,dx,
$$
where $0=\xi_0>\xi_1>\xi_2>\ldots$ satisfy $-p-1<\xi_p<-p$ for
$p=1,2,\ldots$ and $\alpha_p>0, p=0,1,\ldots.$ More precisely, it is
proved that $\xi_p$ is the unique solution $x\in \left]-p-1,-p\right[$ of the
equation $\Psi(1+x)=-\gamma$, and $\alpha_p=1/\Psi'(1+\xi_p)$. Here
$\Psi(x)=\Gamma'(x)/\Gamma(x)$ and $\gamma$ is Euler's constant.

In \cite{B:D2} it is proved directly that $T^n(\underline{0})$
converges to the fixed point $(m_n)$, and this was used to derive that
the same holds independent of where in $\mathcal H$ the iteration
starts.

We cannot apply Banach's fixed point theorem directly because of the
following Lemma.

\begin{lemma}\label{thm:Lipschitz} The best Lipschitz constant $c$ in
$$
d(T(a_n),T(b_n))\le c\,d((a_n),(b_n))\quad\mbox{for}\quad
(a_n),(b_n)\in\mathcal K
$$
is $c=2$.
\end{lemma}
{\it Proof}. For $(a_n),(b_n)$ we find
\begin{eqnarray}\label{eq:basic}
\lefteqn{d(T(a_n),T(b_n))=}\nonumber\\
&&\sum_{k=1}^\infty2^{-k}\frac{\left|\sum_{j=1}^k(b_j-a_j)\right|}{(1+a_1+\ldots+a_k)(1+b_1+\ldots+b_k)}\nonumber\\
&&\le \sum_{k=1}^\infty 2^{-k}\sum_{j=1}^k|a_j-b_j|=\sum_{j=1}^\infty
|a_j-b_j|\sum_{k=j}^\infty 2^{-k}=2d((a_n),(b_n)),
\end{eqnarray}
which shows that $T$ is Lipschitz with constant $c=2$.

Assume next that $T$ satisfies a Lipschitz condition with constant
$c$.

We note that 
$$
d(\underline{0},\underline{a})=\sum_{n=1}^\infty
(a/2)^n=\frac{a}{2-a}.
$$
Furthermore, $T(\underline{0})=\underline{1}$ and for $0\le a<1$ we
have $T(\underline{a})_n=(1-a)/(1-a^{n+1})$, so finally
$$
d(T(\underline{0}),T(\underline{a}))=a\sum_{n=1}^\infty2^{-n}\frac{1-a^n}{1-a^{n+1}}.
$$
This gives
$$
\sum_{n=1}^\infty2^{-n}\frac{1-a^n}{1-a^{n+1}}\le c\frac{1}{2-a},\quad
0<a\le 1,
$$
and letting $a\to 0$ we find  $c\ge 2$. $\quad\square$

\medskip
{\it Proof of Theorem~\ref{thm:attrac}.}

Let us now introduce the set
$$
\mathcal C=\left\{(a_n)\in \mathcal K\mid a_1\ge\tfrac12\right\},
$$
which is a compact convex subset of $\mathcal K$. We first note that
$T(\mathcal K)\subseteq \mathcal C$  because for any $(a_n)\in
\mathcal K$ we have $a_1\le 1$,
hence
$$
T(a_n)_1=\frac1{1+a_1}\ge \frac12.
$$
By \eqref{eq:basic} we always have
\begin{eqnarray}\label{eq:basic2}
\lefteqn{d(T(a_n),T(b_n))=}\nonumber\\
&&\sum_{k=1}^\infty2^{-k}\frac{\left|\sum_{j=1}^k(b_j-a_j)\right|}{(1+a_1+\ldots+a_k)(1+b_1+\ldots+b_k)}\nonumber\\
&&\le\frac1{(1+a_1)(1+b_1)}\sum_{k=1}^\infty 2^{-k}\sum_{j=1}^k|a_j-b_j|=\frac1{(1+a_1)(1+b_1)}\sum_{j=1}^\infty
|a_j-b_j|\sum_{k=j}^\infty 2^{-k}\nonumber\\
&&=\frac2{(1+a_1)(1+b_1)}d((a_n),(b_n)).
\end{eqnarray}
If now $(a_n),(b_n)\in \mathcal C$, we have $a_1,b_1\ge \tfrac12$, and hence
$\frac2{(1+a_1)(1+b_1)}\le\tfrac89$ so that
$$
d(T(a_n),T(b_n))\le\frac{8}{9}d((a_n),(b_n))\quad\mbox{for}\quad
(a_n),(b_n)\in \mathcal C,
$$
showing that $T$ is a contraction on $\mathcal C$. Since $T$ maps
$\mathcal K$ into $\mathcal C$
any fixed point of $T$ on $\mathcal K$ must belong to $\mathcal C$,
and $T:\mathcal C\to \mathcal C$ has a
unique fixed point by Banach's fixed point theorem, and this must be
the sequence $(m_n)$ determined by \eqref{eq:fixrec}. We  also see that
for any $\xi=(a_n)\in \mathcal K$ the iterates
$T^{n+1}(\xi)=T^n(T(\xi))$ converge to the fixed
point $(m_n)$, which finishes the proof of Theorem~\ref{thm:attrac}. 
$\quad\square$

\begin{prop}\label{thm:diam} The compact set $\mathcal H$ defined in 
\eqref{eq:H} has diameter $\diam(\mathcal H)=1$ and the only two
points $a,b\in\mathcal H$ for which $d(a,b)=1$ are 
$\{a,b\}=\{\underline{0},\underline{1}\}$.
\end{prop}

{\it Proof.} It is easy to see that 
$$
1=d(\underline{0},\underline{1})\le \diam(\mathcal H)\le 1,
$$
so we only have to prove that if $d(a,b)=1$ for some $a=(a_n),b=(b_n)$
from $\mathcal H$, then $\{a,b\}=\{\underline{0},\underline{1}\}$. 
Supposing that 
$$
\sum_{n=1}^\infty 2^{-n}|a_n-b_n|=1,
$$
then necessarily $|a_n-b_n|=1$ for all $n\ge 1$. In particular, for
each $n$ necessarily $a_n$ and $b_n$ are either 0 and 1 or these numbers
reversed. Assume now that $a_1=0,b_1=1$. Since $a_n=\int_0^1
x^n\,d\mu(x)$ for a probability measure $\mu$, the condition $a_1=0$
forces $\mu$ to be the Dirac measure $\delta_0$ with unit mass
concentrated at 0, and hence $a_n=0$ for all $n\ge 1$. This shows that
$a=\underline{0}$ and $b=\underline{1}$. $\quad\square$

\noindent Christian Berg\\
Institute of Mathematical Sciences\\
University of Copenhagen\\
Universitetsparken 5\\
DK-2100 K\o benhavn \O, Denmark\\
{\footnotesize E-mail berg@math.ku.dk}

\medskip
\noindent Maryam Beygmohammadi\\
Department of Mathematics\\
Islamic Azad University-Kermanshah branch\\
Kermanshah, Iran\\
{\footnotesize  E-mail maryambmohamadi@mathdep.iust.ac.ir} 
 
\end{document}